\documentclass[11pt]{article}
\usepackage[margin=1in]{geometry}
\usepackage[T1]{fontenc}
\usepackage[utf8]{inputenc}
\usepackage{lmodern,microtype}
\usepackage{amssymb,amsmath,amsthm,graphicx,hyperref,tikz}
\usetikzlibrary{arrows.meta}
\hypersetup{hidelinks}
\newtheorem{theorem}{Theorem}

\newtheorem{lemma}{Lemma}
\newtheorem{claim}{Claim}
\newcommand{\pairs}{\mathcal A}
\definecolor{drawingred}{RGB}{180,35,45}
\definecolor{drawingblue}{RGB}{0,112,163}
\title{At most how many edges can a quasiplanar graph have?}
\author{J\'anos Pach\thanks{Alfréd Rényi Institute of Mathematics, Budapest 1053, Hungary. Supported by ERC Advanced Grant ``GeoScape'', no. 882971. Email: {\tt pach@renyi.hu}.} \and G\'abor Tardos\thanks{Alfréd Rényi Institute of Mathematics, Budapest 1053, Hungary. Supported by ERC Advanced Grant ``ERMiD'', no. 101054936. Email: {\tt tardos@renyi.hu}.} \and Ji Zeng\thanks{Alfréd Rényi Institute of Mathematics, Budapest 1053, Hungary. Supported by ERC Advanced Grants ``GeoScape'', no. 882971 and ``ERMiD'', no. 101054936. Email: {\tt zeng.ji@renyi.hu}.}}
\date{}
\begin{document}
\maketitle

\begin{abstract}
A $k$-quasiplanar graph is a graph drawn in the plane with no $k$ pairwise crossing edges. We construct $k$-quasiplanar graphs on $n$ vertices with $\Omega(k^{1+c}n)$ edges for an absolute constant $c>0$, disproving the widely held belief that the right order of magnitude is $\Theta(kn)$. We prove the lower bound with $c=1/7$, for every integer $k\ge2$ and $n\ge k^{8/7}$. ChatGPT optimized the construction to give $c=323/1827$. We also give a self-contained combinatorial proof of the upper bound $O(k^2n)$, simplifying an earlier proof found by Sadhu using ChatGPT.
\end{abstract}

\section{The results}

A \emph{topological graph} is a simple graph drawn in the plane with distinct points as vertices and Jordan arcs as edges. Edge interiors avoid vertices but may cross each other. We say that two edges \emph{intersect} if the curves have a common interior point. For a positive integer $k$, we say that a topological graph $G$ is \emph{$k$-quasiplanar} if no $k$ distinct edges in $G$ are pairwise intersecting.

It was widely believed that the maximum number of edges in a $k$-quasiplanar graph on $n$ vertices has order $\Theta(kn)$, in the range $2\le k\le n$. Our main result disproves this belief: the dependence on $k$ must be superlinear.

\begin{theorem}\label{lower}
For every integer $k\ge2$ and every integer $n\ge k^{8/7}$, there exists a $k$-quasiplanar topological graph on $n$ vertices with $\Omega(k^{8/7}n)$ edges, where the implied constant is absolute.
\end{theorem}

The recursive construction in Section~\ref{construction} was found with the assistance of ChatGPT. Later, we used ChatGPT to optimize it to give the result with $c=323/1827$. For the sake of simplicity, here we only present the construction and proof for $c=1/7$. The optimal dependence on $k$ remains open.

Pach, Shahrokhi, and Szegedy~\cite{PSSz} conjectured that every $k$-quasiplanar graph on $n$ vertices has $O_k(n)$ edges. Before the recent AI-assisted proofs, this was known only for $k\le4$; see~\cite{AgAP97, Ackerman, AcT07, PRT06}. For fixed $k\ge5$, Fox, Pach, and Suk~\cite{FPS24} proved the upper bound $O_k\bigl(n(\log n)^{2\lceil\log_2 k\rceil-4}\bigr)$.

We also present a self-contained combinatorial proof of the following upper bound.
\begin{theorem}\label{main}
A $k$-quasiplanar topological graph on $n$ vertices has $O(k^2n)$ edges, where the implied constant is absolute.
\end{theorem}

A more complicated proof of this quadratic bound was first found by Sadhu using ChatGPT~\cite{Sadhu}. Our present proof replaces Sadhu's geometric arguments by combinatorial ones. While finalizing this manuscript, we were kindly informed by Istv\'an Tomon that a very similar variant was found by him, jointly with Milojevi\'c, Morawski, and Sudakov~\cite{MMST}. The initial AI-assisted proofs of the $O_k(n)$ bound relied on a theorem of Beeker and Lazarovich~\cite{BL}, whose proof involves advanced topological techniques and combinatorial arguments. The proof given here avoids those techniques.

\bigskip \noindent{\bf Acknowledgements.} This paper grew out of an effort to interpret and simplify a proof generated by ChatGPT-6. While preparing~\cite{PTZ}, we discovered Sadhu's preprint~\cite{SadhuV1} offering a very similar strategy. Both contained combinatorial arguments inspired by~\cite{BL}, which led only to an $O_k(n)$ upper bound. The lower-bound construction and its optimization were also developed with the assistance of ChatGPT.

\section{Set-up for the upper bound}

Let $G$ be a $k$-quasiplanar topological graph with $n$ vertices,
and put $d=k-1$. The cases $k=1$ and $n\le2$ are immediate, so assume $k\ge2$ and $n\ge3$.
Without loss of generality, we assume that the edges are not
self-intersecting and consist of finitely many straight line segments.
We may also assume that all intersections of edge interiors are proper
crossings and that no three edge interiors meet at one point; these
conditions can be achieved without introducing new intersecting pairs.
We choose a triangular region $T$ in the plane that contains all
vertices and edges of $G$ in its interior. For each vertex $u$ of $G$,
we choose a small triangle $T_u$ containing $u$; we call these
\emph{special triangles}. Their interiors will be called \emph{holes}.

We triangulate $T$ using its vertices and those of the special
triangles, without introducing any additional vertices, so that
every special triangle is a face. We call the remaining faces
\emph{regular triangles}. We call the edges of the triangulation
\emph{triangulation edges} and the edges of $G$ \emph{drawing edges}.
By Euler's formula, we obtain $O(n)$ regular triangles and $O(n)$
triangulation edges. More precisely, the number of regular triangles is
\begin{equation}\label{triangles}
                         t=5n+1.
\end{equation}
By choosing the $T_u$ sufficiently small and slightly perturbing $G$,
we may assume that the special triangles are pairwise disjoint,
contain no crossings, and satisfy the following conditions:
\begin{enumerate}
\item No drawing edge contains a triangulation vertex, and no
triangulation edge contains a crossing point between drawing edges.
\item Each intersection between a drawing edge and a triangulation
edge is a proper crossing.
\item The drawing edge $uv$ crosses each of $\partial T_u$ and
$\partial T_v$ exactly once and meets no other special triangle.
\end{enumerate}

\begin{figure}[tb]
\centering
\begin{tikzpicture}[x=0.95cm,y=0.95cm,line cap=round,line join=round,
  every node/.style={font=\small}]
  \begin{scope}[shift={(0,0)}]
    \fill[black!5] (-1.50,0)
      .. controls (-0.50,1.60) and (0.50,1.60) .. (1.50,0)
      -- cycle;
    \draw[black,line width=0.8pt] (-2.65,0)--(2.65,0);
    \draw[drawingred,line width=1.25pt] (-2,{-14/15})
      .. controls ({-2/3},{86/45}) and ({2/3},{86/45}) .. (2,{-14/15});
    \draw[drawingblue,line width=1.8pt] (-0.85,-0.70)
      .. controls (-0.35,{29/30}) and ({2/3},1.80) .. (2.20,1.80);
    \fill (-1.50,0) circle (1.1pt) (1.50,0) circle (1.1pt);
    \node[above] at (-0.65,1.10) {$\alpha$};
    \node[drawingblue,below] at (1.30,1.53) {$\beta$};
    \node at (0.35,0.40) {$L$};
  \end{scope}
  \draw[-{Stealth[length=2.5mm]},line width=0.9pt]
    (3.00,0.43)--(3.85,0.43);
  \begin{scope}[shift={(6.85,0)}]
    \draw[black,line width=0.8pt] (-2.65,0)--(2.65,0);
    \draw[drawingred,line width=1.25pt] (-2,{-14/15})
      .. controls ({-11/6},{-4/9}) and ({-5/3},-0.20) .. (-1.50,-0.20)
      -- (1.50,-0.20)
      .. controls ({5/3},-0.20) and ({11/6},{-4/9}) .. (2,{-14/15});
    \draw[drawingblue,line width=1.8pt] (-0.85,-0.70)
      .. controls (-0.35,{29/30}) and ({2/3},1.80) .. (2.20,1.80);
    \node[drawingblue,below] at (1.30,1.53) {$\beta$};
    \node[below] at (0.35,-0.24) {$\alpha'$};
  \end{scope}
\end{tikzpicture}
\caption{We can eliminate the lune $L$ by replacing $\alpha$ by $\alpha'$. A drawing edge $\beta$ that crosses $\alpha'$ must also cross $\alpha$ by minimality of the lune $L$.}
\label{lune}
\end{figure}
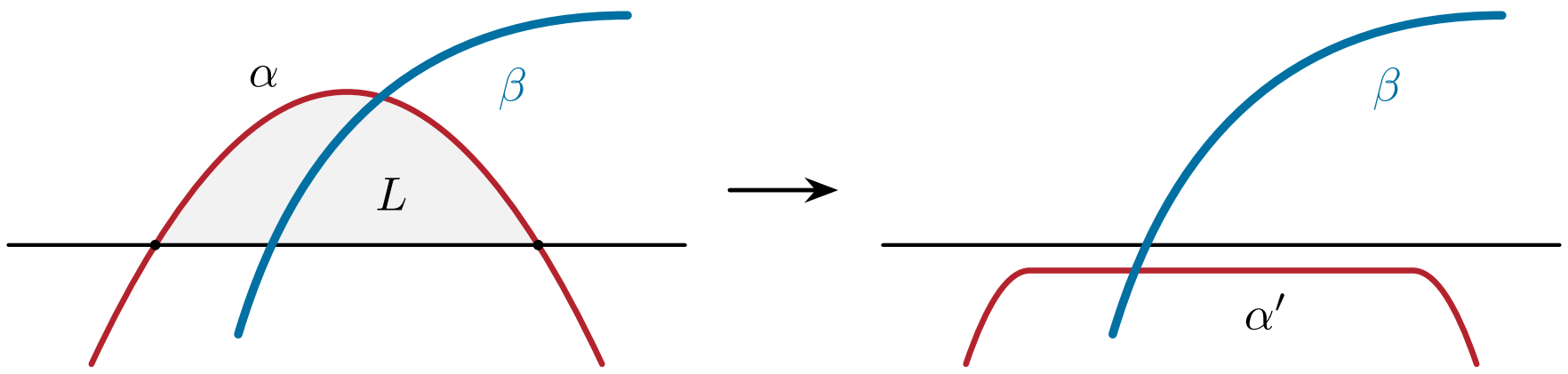

We can further assume without loss of generality that
\begin{enumerate}
\setcounter{enumi}{3}
\item\label{normal} Any two consecutive intersection points with
triangulation edges along a drawing edge $uv$ lie on different
triangulation edges.
\end{enumerate}
Indeed, if this condition fails, some drawing edge $uv$ crosses a
triangulation edge $pq$ at two points bounding an arc $\alpha$ that
crosses no other triangulation edge. Then $pq$ and $\alpha$ enclose
a bounded region in the plane, called a \emph{lune}. We choose a
lune minimal under containment among all such lunes, and redraw the
arc $\alpha$ on $uv$ as shown in Figure~\ref{lune}. This eliminates
the lune and, by its minimality, introduces no new pairs of crossing
drawing edges. The same minimality excludes other portions of $uv$
from the lune, so the redrawn edge remains simple. Each move is
made outside the special triangles and removes two intersections
with triangulation edges. Hence we can repeatedly redraw $G$ to
guarantee Condition~\ref{normal} while preserving $k$-quasiplanarity.

For an oriented triangulation edge $pq$, let $X_{pq}$ be the set
of pairs $(a,e)$ with $e$ a drawing edge and $a\in e\cap pq$.
We call these pairs \emph{occurrences}; different intersections of
the same drawing edge give different occurrences. We define two
linear orders $L_{pq}$ and $R_{pq}$ on $X_{pq}$ as follows.

For $(a,e)\in X_{pq}$ follow the drawing edge $e$ from $a$ starting at a fixed side (left or right) of $pq$. Whenever $e$ enters a regular triangle $xyz$ through the directed
side $xy$, record $0$ if it exits through $xz$, and $1$ if it exits
through $zy$. Condition~\ref{normal} guarantees these are the only
two possibilities. Give the next side the orientation
$x\to z$ or $z\to y$, respectively. Continue this process until $e$
reaches a hole boundary after finitely many crossings. Also record this terminal point at the hole boundary. This records a finite binary word and a terminal point for each element of $X_{pq}$. If the chosen side initially
faces a hole, the word is empty and the terminal point is $a$ itself.

Order these words lexicographically, with $0<1$. Note that $X_{pq}$ may have several elements with the same word, but the terminal points for these elements are all on the same directed edge. To order these elements, use the order of their terminal points along that directed edge.
This defines $L_{pq}$ when the chosen side is the left side of $pq$,
and $R_{pq}$ when it is the right side. On a hole side, one of these
is just the physical order of the occurrences from $p$ to $q$.
Outer boundary sides have no occurrences.

For a fixed starting side and choice of left or right, no word is a
proper prefix of another: once a sequence of crossings reaches a hole,
every drawing edge following that sequence terminates there. A word and its terminal point
determine the element of $X_{pq}$, so
the tie-breaking rule gives genuine linear orders on $X_{pq}$.

A subset of a linearly ordered set is a \emph{prefix} if it contains
every element preceding any of its members, and a \emph{suffix} if
it contains every element following any of its members. The empty
set and the whole ordered set are both prefixes and suffixes.

\begin{claim}\label{orders}
If the orders of $(a,e),(b,f)\in X_{pq}$ differ in $L_{pq}$ and $R_{pq}$, then the drawing edges $e$ and $f$ cross. Therefore, no sequence of $d+1$ occurrences is increasing in $L_{pq}$ and
decreasing in $R_{pq}$.

For a regular triangle $xyz$, consider the map $(a,e)\mapsto(b,e)$ whenever $a\in xy\cap e$, $b\in xz\cap e$ and the segment of $e$ between $a$ and $b$ is inside the triangle $xyz$. This map is a bijection between subsets of $X_{xy}$ and $X_{xz}$; it preserves both linear orders $L_{\cdot}$ and $R_{\cdot}$ and the subsets are prefixes in one of these orders. Similarly, the map $(a,e)\mapsto(b,e)$ if $a\in xy\cap e$, $b\in zy\cap e$, and the segment of $e$ between $a$ and $b$ is inside the triangle $xyz$ is a bijection between subsets of $X_{xy}$ and $X_{zy}$ preserving both linear orders. In this case the subsets are suffixes in one of the two orders.
\end{claim}

\begin{figure}[htbp]
  \centering
  \begin{tikzpicture}[scale=0.8,
    x=1cm,y=1cm,
    line cap=round,line join=round,
    every node/.style={font=\normalsize},
    strand e/.style={
      draw=cyan!65!blue,
      line width=1.2pt,
      dash pattern=on 3pt off 2.5pt
    },
    strand f/.style={
      draw=red!80!black,
      line width=1.8pt
    }
  ]
    \coordinate (p) at (-3,0);
    \coordinate (q) at ( 3,0);
    \coordinate (r) at ( 0,4.5);
    \coordinate (a) at (-0.9,0);
    \coordinate (b) at ( 0.8,0);

    \draw[black,line width=0.8pt]
      (p)--(q)--(r)--cycle;

    \draw[strand e] (-0.9,-0.55)--(a)
      .. controls (-0.9,0.70) and (-1.55,1.45)
      .. (-2.80,1.65);

    \draw[strand f] (0.8,-0.55)--(b)
      .. controls (0.8,1.10) and (-0.25,2.35)
      .. (-1.90,3.00);

    \fill[cyan!65!blue] (a) circle (1.3pt);
    \fill[red!80!black] (b) circle (1.3pt);

    \node[below left,cyan!65!blue] at (a) {$a$};
    \node[below right,red!80!black] at (b) {$b$};
    \node[cyan!65!blue] at (-1.10,1.02) {$e$};
    \node[red!80!black] at (0.48,1.70) {$f$};

    \node[below left]  at (p) {$p$};
    \node[below right] at (q) {$q$};
    \node[above]       at (r) {$r$};
  \end{tikzpicture}

  \caption{Since $e$ and $f$ do not cross, $e$ must leave the triangle through $pr$.}
  \label{downset}
\end{figure}
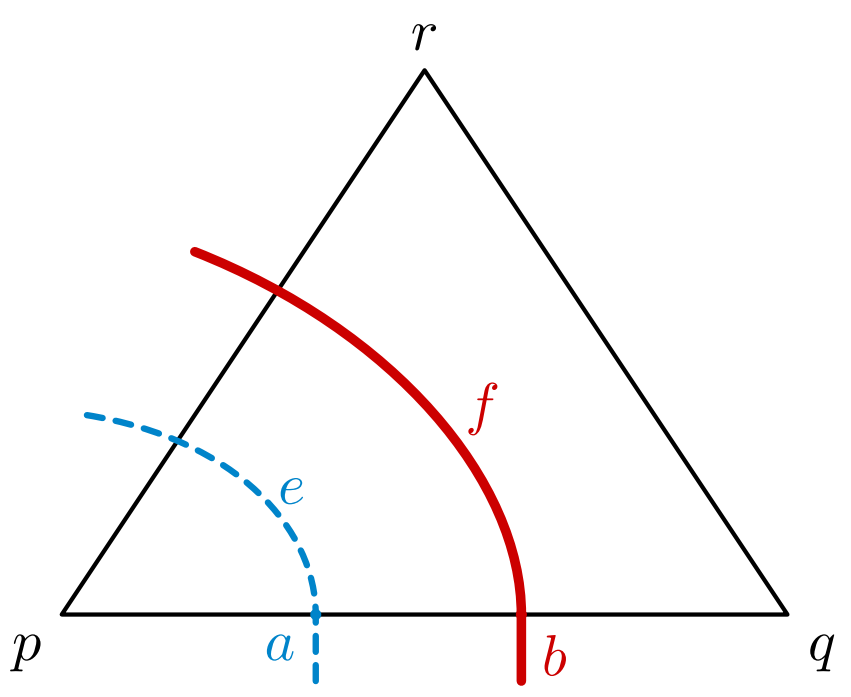

\begin{proof}
In a triangle $xyz$, disjoint drawing edges entering through $xy$ meet the
complementary boundary path $xz\cup zy$ in the same $x$-to-$y$ order;
see Figure~\ref{downset}.
Apply this observation to two occurrences belonging to drawing edges
which do not cross, or to the same drawing edge. While their words
agree, their continued components are distinct and disjoint.
Indeed, identical components would give identical starting occurrences
by backward continuation. At the first different exit, the $0$-arc
precedes the $1$-arc; if the words are equal, their terminal points
retain the original order. Hence, for $(a,e),(b,f)\in X_{pq}$ with $e=f$ or $e$ and $f$ have disjoint interiors, their orders in both $L_{pq}$ and $R_{pq}$ agree with the order of $a$ and $b$ on $pq$. This proves the first claim. The claim on monotone sequences follows because a pairwise crossing family of drawing edges has at most $d$ elements.

Suppose without loss of generality that $xyz$ lies to the left of $xy$. The first map maps the elements of $X_{xy}$ whose word for $L_{xy}$ starts with a $0$; these elements form a prefix in $L_{xy}$. The second map maps the elements of $X_{xy}$ whose word for $L_{xy}$ starts with a $1$; these elements form a suffix in $L_{xy}$. If either map maps $(a,e)$ to $(b,e)$, the word of $(b,e)$ for $L_{xz}$ or $L_{zy}$ is obtained from the word for $(a,e)$ by deleting this initial $0$ or $1$, while the terminal point remains the same. This proves that the maps preserve the order $L_{\cdot}$. The inverse of the first map maps the elements $(b,e)\in X_{xz}$ whose words for $R_{xz}$ start with a $0$; these elements form a prefix in $R_{xz}$. Similarly, the inverse of the second map maps $(b,e)\in X_{zy}$ whose word for $R_{zy}$ starts with a $1$; these elements form a suffix in $R_{zy}$. The word for $R_{xy}$ of the preimage is obtained by dropping this initial digit and the terminal point does not change. This shows that both maps preserve the orders $R_{\cdot}$ too.
\end{proof}

\section{Two lemmas on ordered sets}\label{combinatorics}

Let $X$ be a finite set with two linear orders $L,R$, and suppose
that no sequence of $d+1$ elements is increasing in $L$ and decreasing
in $R$. For $x\in X$, define
\[
 h_X(x)=\max\{s:x_1<_L\cdots<_L x_s=x,
                    \ x_1>_R\cdots>_R x_s\},
 \qquad C_i(X)=\{x:h_X(x)=i\}.
\]
Order each nonempty $C_i(X)$ by $L$ and call it a \emph{rank chain}.
If $x<_Ly$ and $y<_Rx$, then $h_X(y)>h_X(x)$ and therefore $x$ and $y$ cannot be in the same rank chain. This shows that $L$ and $R$ agree on each chain. Write
$\pairs(X)$ for the pairs $(x,y)$ such that $y$ immediately follows
$x$ in one of these chains. For $Y\subseteq X$, give $Y$ the induced
orders and recompute its ranks; write
$\pairs(X)|_Y=\pairs(X)\cap(Y\times Y)$.

\begin{lemma}\label{chains}
The rank chains partition $X$ into at most $d$ chains increasing in
both orders. At most $d$ elements $x\in X$ has no $y$ with $(x,y)\in\pairs(X)$.
If $Y$ is a prefix or suffix of either order, at most $d$
pairs in $\pairs(X)$ have exactly one member in $Y$.
\end{lemma}

\begin{proof}
If $x<_L y$ and $x>_R y$, then $h_X(y)\ge h_X(x)+1$.
Thus $<_L$ and $<_R$ agree on each rank chains. There are at
most $d$ such chains as we assumed that no $d+1$ elements are increasing in $L$ and decreasing in $R$, so $h_X$ is bounded by $d$. Only the last element $x$ of a rank chain has no $y$ with $(x,y)\in\pairs(X)$. Finally, a prefix or suffix of either
order meets each chain in a prefix or suffix, cutting at
most one consecutive pair.
\end{proof}

\begin{claim}\label{interval}
    Suppose $Y$ is a prefix or suffix of $X$ under either $L$ or $R$. Then $C_i(X) \cap C_j(Y)$ is an interval in both chains $C_i(X)$ and $C_j(Y)$, where the ranks on $Y$ are computed using only elements of $Y$.
\end{claim}

\begin{proof}
    First, it suffices to consider the cases where $Y$ is a prefix or suffix under $L$. Once these cases are proved, we can replace the orders $(L,R)$ by $(R^{\mathrm{op}},L^{\mathrm{op}})$, where $\mathrm{op}$ denotes the reversed order, to conclude the other cases.

    Suppose that $Y$ is a prefix of $X$ under $L$. For $x \in Y$, a sequence $x_1<_L\cdots<_L x_s=x$ with $x_1>_R\cdots>_R x_s$ and $s = h_X(x)$ lies entirely in $Y$. Thus $h_X(x) = h_Y(x)$ for $x\in Y$, which implies $C_i(Y) = C_i(X) \cap Y$. Hence $C_i(X) \cap C_j(Y)$ equals $C_i(Y)$ when $i = j$ and $\emptyset$ otherwise. The set $C_i(Y)$ is an interval in itself, and $C_i(X) \cap Y$ is a prefix of $C_i(X)$. This proves the prefix case.

    Now suppose that $Y$ is a suffix of $X$ under $L$. For $S = X$ or $S = Y$, we describe an algorithm to compute $h_S(v)$ for all $v\in S$. We shall then compare the algorithm on $X$ and $Y$ to prove the claim. We process $v\in S$ in increasing order under $L$, maintaining a subset $H\subseteq S$ called \emph{state}, which we list in decreasing order under $R$. Initially, $H = \emptyset$. Before processing $v\in S$, write $H$ as $t_1 >_R \ldots >_R t_r$ and let $j = 1+\max\bigl(\{i:~ t_i >_R v\}\cup\{0\}\bigr)$. If $j\leq r$, we replace $t_j$ with $v$; if $j = r+1$, we add $v$ as $t_j$ to $H$. Inductively, before each step, $t_i$ is the $R$-largest processed element of rank $i$. Thus $h_S(v) = j$, since some $L$-increasing, $R$-decreasing sequence of length $i$ can be extended to $v$ exactly when $t_i>_R v$.

    Let $H_S(v)$ be the state $H$ immediately after $v$ is processed. We first process the elements of $X\setminus Y$ in the $X$-process, while the $Y$-process has empty state. We then process the elements of $Y$ in both processes (recall $Y$ is a suffix). We check that, as the algorithms proceed, we maintain $H_Y(v) \subseteq H_X(v)$ and \begin{equation}\label{algorithm}
        \text{elements in $H_X(v)\setminus H_Y(v)$ are either removed or replaced by smaller elements under $R$.}
    \end{equation} Indeed, the inclusion holds initially because the $Y$-state is empty. Before processing $v\in Y$, write the two states as $A,B$ for $X,Y$ respectively, and let $a = \max_R\{x\in A:~x<_R v\}$ and $b = \max_R\{x\in B:~x<_R v\}$ whenever these maxima exist. Since inductively $B\subseteq A$, we cannot have $a<_R b$, or $b$ existing without $a$. The updated states are $H_X(v)=(A\setminus\{a\})\cup\{v\}$ and $H_Y(v)=(B\setminus\{b\})\cup\{v\}$, with no deletion when the corresponding maximum does not exist. If $a=b$, or neither exists, the difference $A\setminus B$ is unchanged. If $a>_R b$, then $a\notin B$, and the new difference replaces $a$ with $b$. If $a$ exists and $b$ does not, then $a\notin B$, and the new difference removes $a$. These cases preserve the inclusion and verify \eqref{algorithm}.

    Since $h_S(v) = |\{x \in H_S(v):~ x \geq_R v\}|$ for $S=X$ or $Y$, and $H_Y(v) \subseteq H_X(v)$, we have \begin{equation*}
        h_X(v)-h_Y(v) = |\{x \in H_X(v)\setminus H_Y(v):~ x \geq_R v\}|
    \end{equation*} Together with \eqref{algorithm}, this gives \begin{equation*}
        h_X(v)-h_Y(v) \leq |\{x \in H_X(u)\setminus H_Y(u):~ x \geq_R v\}| \leq h_X(u)-h_Y(u)
    \end{equation*} if $u <_L v$ and $u <_R v$ simultaneously for $u,v\in Y$. Now consider $h_Y(v)$ restricted to $C_i(X) \cap Y$, where the two orders $L$ and $R$ coincide. The inequality above implies that it is nondecreasing, since $h_X(v)$ is fixed. Hence, $C_i(X)\cap C_j(Y)$ is an interval in $C_i(X)\cap Y$, and since $Y$ is a suffix, it is also an interval in $C_i(X)$. Similarly, $h_X(v)$ is nonincreasing on $C_j(Y)$, so $C_i(X)\cap C_j(Y)$ is an interval in $C_j(Y)$.
\end{proof}

\begin{lemma}\label{transfer}
Let $X_1,X_2$ be finite sets with linear orders $L_i,R_i$. Assume that neither set has a sequence of $d+1$ elements forming an increasing sequence under one a decreasing sequence under the other order.
For $i=1,2$, let $Y_i\subseteq X_i$ be a prefix or suffix of
one of the two orders, and let $\phi:Y_1\to Y_2$ be a bijection preserving both
orders. Then at most $O(d^2)$ pairs of $\pairs(X_1)|_{Y_1}$ have images
outside $\pairs(X_2)|_{Y_2}$. The same bound holds in the reverse
direction.
\end{lemma}

\begin{proof}
    Each pair in $\pairs(X_1)|_{Y_1}$ but not in $\pairs(Y_1)$ joins different intervals $C_i(X_1)\cap C_j(Y_1)$ from Claim~\ref{interval}. There are at most $d^2$ such nonempty intervals, each with at most one outgoing consecutive pair, so there are $O(d^2)$ such pairs. $\phi$ identifies $\pairs(Y_1)$ with $\pairs(Y_2)$. The number of pairs in $\pairs(Y_2)$ but not in $\pairs(X_2)|_{Y_2}$ is at most $O(d^2)$ by the same argument. This proves the lemma.
\end{proof}

\section{Proof of the upper bound}

Use the two orders of Claim~\ref{orders} to form the rank chains on
every oriented triangulation edge $pq$. If one occurrence immediately
follows another in its rank chain, call it its \emph{successor}.
A \emph{representative} of an edge $e$ of $G$ is an occurrence
$(a,e)\in X_{pq}$, on an oriented side $pq$ of a regular triangle
$pqr$, such that one of the following holds:
\begin{description}
\item[(a)] $(a,e)$ has no successor in its rank chain, or
\item[(b)] $(a,e)$ has a successor $(b,f)$, and $e,f$ leave $pqr$
through different sides, or
\item[(c1)] $(a,e)$ has a successor $(b,f)$, both leave $pqr$ through
$pr$, but their images, in this order, are not consecutive in a
rank chain on the oriented side $pr$, or
\item[(c2)] $(a,e)$ has a successor $(b,f)$, both leave $pqr$ through
$rq$, but their images, in this order, are not consecutive in a
rank chain on the oriented side $rq$.
\end{description}

We argue that every drawing edge $e=uv$ has a representative.
Let $a_0$ be its unique crossing with $\partial T_u$, on an oriented
triangulation edge $p_0q_0$. If $(a_0,uv)$ has no successor, it is
a representative of type (a). Otherwise, let $(b_0,uw)$ be its
successor. Every occurrence on this hole side belongs to an edge
incident with $u$, and $v\ne w$ because $G$ is simple.

We construct a sequence of pairs of occurrences
\[
                   \bigl((a_i,uv),(b_i,uw)\bigr)
                   \quad\text{on }p_iq_i
\]
by following $uv$ and $uw$ through the triangles. More precisely,
let $p_iq_ir_i$ be the regular triangle which the two drawing edges,
oriented from $u$, enter at $a_i,b_i$. If both leave through $p_ir_i$,
let $a_{i+1},b_{i+1}$ be their next crossings and put
$p_{i+1}q_{i+1}=p_ir_i$. If both leave through $r_iq_i$, define the
next terms similarly, with $p_{i+1}q_{i+1}=r_iq_i$. These are the
order-preserving bijections of Claim~\ref{orders}.

The construction ends when the two edges take different exits.
Such an ending exists because both routes have finitely many
crossings and end in different special triangles. Let $s$ be this
last index. If $(a_s,uv)$ and $(b_s,uw)$ are consecutive in a rank
chain, in this order, then $(a_s,uv)$ is a representative of type (b).
Otherwise, since $(b_0,uw)$ is a successor of $(a_0,uv)$, there is some
$0\le i<s$ for which the pair is in $\pairs(X_{p_iq_i})$ but its image is not in $\pairs(X_{p_{i+1}q_{i+1}})$. Then $(a_i,uv)$ is
a representative of type (c1) or (c2).

Finally, we bound the number of representatives. Fix an oriented
side $pq$ of a regular triangle $pqr$. Since $X_{pq}$ has at most
$d$ rank chains by Lemma~\ref{chains}, there are at most $d$ representatives of type (a).
The two exit families are a prefix and a suffix of the
order facing $pqr$, so Lemma~\ref{chains} gives at most $d$
representatives of type (b). By Claim~\ref{orders} and
Lemma~\ref{transfer}, there are $O(d^2)$ representatives of
each of types (c1) and (c2).

There are six oriented side incidences per regular triangle, and
$t=5n+1$ regular triangles. Every drawing edge has a representative,
and a representative identifies its drawing edge. Therefore
\[
                         |E(G)|=O(td^2)=O(k^2n).
\]
This proves Theorem~\ref{main}.

\section{A superlinear lower bound}\label{construction}

We prove Theorem~\ref{lower} by a recursive construction of drawings of complete bipartite graphs. The construction yields the exponent $1+c=8/7$ and shows that an $O(kn)$ upper bound is impossible.

We prove Theorem~\ref{lower} by a recursive construction of drawings of complete bipartite graphs. Let $\delta=1/4$. We choose this parameter for simplicity. For positive integers $a$ and $z$, we define the topological graph $G_{a,z}$, a drawing of the complete bipartite graph $K_{a,z}$. We will maintain the following properties.
\begin{enumerate}
\item $G_{a,z}$ is contained in an axis-aligned rectangle $R$,
\item the $a$ left vertices of $G_{a,z}$ are in the interior of $R$, \item the $z$ right vertices are on the right edge of $R$,
\item the edges are $x$-monotone and consist of finitely many straight line segments,
\item two edges incident to the same right vertex do not cross.
\end{enumerate}

The construction is recursive, starting with the case $a\le3$, where we place the left vertices on a vertical line in the interior of $R$ and join them by straight line segments to the right vertices. These drawings satisfy conditions (1)--(5), and every family of pairwise crossing edges has size at most $a\le3$.

For $a\ge4$ we set $b=\lfloor\delta a\rfloor$ and $c=a-b$. We start with the topological graphs $G_{b,z}$ and $G_{c,z}$ and place them in the plane such that the enclosing rectangle $R'$ of the former is to the left of the enclosing rectangle $R''$ of the latter. We select an enclosing rectangle $R$ for $G_{a,z}$ that contains both $R'$ and $R''$ in its interior. The set of left vertices of $G_{a,z}$ is the union of the left vertices of $G_{b,z}$ and $G_{c,z}$. For the right vertices we select $v_1,\dots,v_z$ on the right edge of $R$ indexed so their $y$-coordinates increase. Let $v'_1,\dots,v'_z$ and $v''_1,\dots,v''_z$ be the right vertices of $G_{b,z}$ and $G_{c,z}$, respectively, also indexed in increasing order of their $y$-coordinates. We will modify each edge of both $G_{b,z}$ and $G_{c,z}$ to obtain an edge of $G_{a,z}$ as follows.

Take a vertical line $l''$ to the left of the right edge of $R''$ but close enough so that all crossings and all left vertices of $G_{c,z}$ lie to the left of it. Now take an edge of $G_{c,z}$ incident to a right vertex $v''_i$ and extend its part left of $l''$ with a single straight line segment going to $v_i$. We call the edges of $G_{a,z}$ obtained this way the \emph{black edges}. Note that two black edges cross each other if and only if the corresponding edges of $G_{c,z}$ cross.

Take a vertical line $l'$ to the left of the right edge of $R'$ but close enough so that all crossings and all left vertices of $G_{b,z}$ lie to the left of it. Take an edge of $G_{b,z}$ incident to a right vertex $v'_i$. If $i\le z/2$, we extend the part of this edge to the left of $l'$ by an $x$-monotone curve passing strictly under $R''$ and ending at $v_{2i}$. We call the edges of $G_{a,z}$ obtained this way the \emph{red edges}. We do this in such a way that the red edges consist of straight line segments that do not leave $R$, two red edges cross each other if and only if the corresponding edges of $G_{b,z}$ cross, while a red edge to $v_i$ crosses a black edge to $v_j$ if and only if $i>j$. Note that the right vertices of red edges have even indices.

Finally, for an edge of $G_{b,z}$ incident to a right vertex $v'_i$ with $i>z/2$, we extend the part of it to the left of $l'$ by an $x$-monotone curve passing strictly over $R''$ and ending at $v_{2i-z-1}$ if $z$ is even or at $v_{2i-z}$ if $z$ is odd. We call the edges of $G_{a,z}$ obtained this way the \emph{blue edges}. We do this in such a way that the blue edges consist of straight line segments that do not leave $R$, two blue edges cross each other if and only if the corresponding edges of $G_{b,z}$ cross, while a blue edge to $v_i$ crosses a black edge to $v_j$ if and only if $i<j$. Note that the right vertices of blue edges have odd indices.

Clearly, conditions (1)--(5) are maintained by this construction. It remains to bound the number of pairwise crossing edges.

Let $\alpha=1/8$, $\beta=3/4$, and let $C$ be a sufficiently large absolute constant.

\begin{lemma}\label{crossing-interval}
Let $a$, $z$, and $h$ be positive integers with $h\le z$. Let $H$ be a set of pairwise crossing edges in $G_{a,z}$ such that their right vertices are within an interval $I$ of length $h$, that is $I=\{v_i\mid s\le i<s+h\}$ for some $1\le s\le z-h+1$. Then
\[
|H|\le Ca^\alpha h^\beta.
\]
\end{lemma}

\begin{proof}
We prove the lemma by induction on $a$. In the base case $a\le3$ we have at most three pairwise crossing edges, so the claim holds if $C\ge3$.

For $a\ge4$ we follow the recursive construction of $G_{a,z}$ from $G_{b,z}$ and $G_{c,z}$. As $b,c<a$ we assume the statement of the lemma holds for them. Let $H$ consist of $t_1$ black, $t_2$ red and $t_3$ blue edges. Let $H_1$ consist of the edges of $G_{c,z}$ corresponding to the black edges in $H$. Similarly, let $H_2$ and $H_3$ consist of the edges of $G_{b,z}$ corresponding to the red and blue edges in $H$, respectively. Note that all three of these families are pairwise crossing. Empty intervals are allowed below, with zero contribution from the corresponding empty family. We consider two cases according to whether $t_1=0$.

Suppose $t_1=0$, that is, $H$ consists of red and blue edges only. The right vertices of the edges in $H_2$ lie in an interval of length $h'$, where $h'$ is the number of even-indexed vertices in $I$. Similarly, the right vertices of $H_3$ come from an interval of length $h''$, the number of odd-indexed vertices in $I$. We have $h'+h''=h$ and the inductive hypothesis gives
\[
|H|=t_2+t_3\le Cb^\alpha h'^\beta+Cb^\alpha h''^\beta\le2C(\delta a)^\alpha(h/2)^\beta=Ca^\alpha h^\beta
\]
as needed. The second inequality follows from $b\le\delta a$ and the concavity of $x^\beta$. The last equality comes from $2\delta^\alpha/2^\beta=1$.

Suppose now that $t_1\ge 1$, that is, there is a black edge in $H$. As all red and black edges cross in $H$, the right vertices of the red edges must all be strictly above the right vertices of the black edges. Similarly the right vertices of blue edges must all be below the right vertices of black edges. Thus, the interval $I$ can be split into three pairwise disjoint intervals $I_1$, $I_2$, and $I_3$ such that the right vertices of black edges are in $I_1$, the right vertices of red edges are in $I_2$ and the right vertices of blue edges are in $I_3$. For the length $h_i$ of the interval $I_i$ we have $h_1+h_2+h_3=h$. The right vertices of $H_1$ are also in an interval of length $h_1$. The right vertices of $H_2$ correspond to even-indexed elements of $I_2$, so they are from an interval of length at most $\lceil h_2/2\rceil$. Similarly, the right vertices of $H_3$ are from an interval of length at most $\lceil h_3/2\rceil$. The inductive hypothesis yields
\[
|H|=t_1+t_2+t_3\le Cc^\alpha h_1^\beta+Cb^\alpha\lceil h_2/2\rceil^\beta+Cb^\alpha\lceil h_3/2\rceil^\beta.
\]
Fix $\eta>1$ with $\sqrt3/2+\eta^3/8<1$.
For $i=2,3$, we have $t_i\le Cb^\alpha(\eta h_i/2)^\beta$:
for $h_i\ge1/(\eta-1)$ use $\lceil h_i/2\rceil\le\eta h_i/2$,
and for bounded $h_i$ use $t_i\le h_i$, since no two edges of $H$ have the same right vertex, and enlarge $C$. This choice is uniform in $a$, since $b\ge1$.
H\"older's inequality gives
\[
 |H|\le C\left(\sqrt c+\frac{\eta^3}{4}\sqrt b\right)^{1/4}
 h^{3/4}<Ca^\alpha h^\beta
\]
for sufficiently large $a$, since $b/a\to1/4$.
To cover bounded $a$, first note that an independent induction bounds
every crossing family by $a^2$: the base cases are immediate, and
the three color classes give the bound $c^2+2b^2\le a^2$.
Thus, if $a_0$ is a fixed threshold above which the strict inequality
holds, choosing $C\ge a_0^{2-\alpha}$ covers every $a\le a_0$
and every $h\ge1$. All the required choices of $C$ can be made
before the induction begins.

This concludes the proof of the lemma.
\end{proof}

\begin{proof}[Proof of Theorem~\ref{lower}]
Given $k>C$, let $n_k$ be the largest positive integer such that $Cn_k^{\alpha+\beta}<k$. Clearly $n_k=\Theta(k^{1/(\alpha+\beta)})$. The topological graph $G_{n_k,n_k}$ is $k$-quasiplanar by Lemma~\ref{crossing-interval} and has $2n_k$ vertices and $n_k^2=\Theta(k^{1/(\alpha+\beta)}n_k)$ edges. For $n\ge2n_k$ one can take the disjoint union of $\lfloor n/(2n_k)\rfloor$ copies of $G_{n_k,n_k}$ and isolated vertices to obtain a $k$-quasiplanar graph on $n$ vertices with $\Theta(k^{1/(\alpha+\beta)}n)$ edges. Since $C\ge2$, we have $2n_k<2(k/C)^{8/7}\le k^{8/7}$, so this covers the claimed range of $n$. The remaining values $2\le k\le C$ are handled by a planar star.
\end{proof}

\end{document}